\documentclass{article}
\usepackage[dvips]{graphicx}
\usepackage{amssymb,amsfonts,latexsym,amsmath,amscd}

\title{On Roots of Eigenpolynomials for Degenerate Exactly-Solvable Differential Operators}
\author{Tanja Bergkvist \& Jan-Erik Bj\"ork\\ Department of Mathematics, University of Stockholm,\\
S-106 91 Stockholm, Sweden\\ e-mail: tanjab@math.su.se \& jeb@math.su.se}
\date{}
\begin{document}
\maketitle
\begin{abstract}
In this paper we partially settle our conjecture from \cite{B} on the roots of eigenpolynomials for degenerate exactly-solvable operators. Namely, for any such operator we establish a lower bound (which supports our conjecture) for the largest modulus of all roots of its unique and monic eigenpolynomial $p_n$ as the degree $n$ tends to infinity. The main theorem below thus extends earlier results obtained in \cite{B} for a restrictive class of operators.
\end{abstract}
\section{Introduction} We are interested in roots of eigenpolynomials satisfying certain linear differential equations. Namely, consider an operator 
$$T=\sum_{j=1}^{k}Q_jD^j$$ where $D=d/dz$ and the $Q_j$ are complex polynomials in one variable satisfying the condition $\deg Q_j\leq j$, with equality for at least one $j$, and in particular $\deg Q_k<k$ for the leading term. Such operators are referred to as \textit{degenerate exactly-solvable operators}\footnote{
Correspondingly, operators for which $\deg Q_k=k$ are called \textit{non-degenerate exactly-solvable operators}. We have treated roots of eigenpolynomials for these operators in \cite{BR}.}, see \cite{B}. We are interested in eigenpolynomials of $T$, that is polynomials satisfying 
\begin{equation}\label{eigenvalueproblem}
T(p_{n})=
\lambda_{n}p_{n}
\end{equation}
 for some value of the spectral parameter $\lambda_{n}$, where $n$ is a positive integer and $\deg p_n=n$. The importance of studying eigenpolynomials for these operators is among other things motivated by numerous examples coming from classical orthogonal polynomials, such as the Laguerre and Hermite polynomials, which appear as solutions to (\ref{eigenvalueproblem}) for certain choices on the polynomials $Q_j$ when $k=2$. Note however that for the operators considered here the sequence of eigenpolynomials $\{p_n\}$ is in general \textit{not} an orthogonal system.\\
Let us briefly recall our previous results:\\\\
\textbf{A.} In \cite{BR} we considered eigenpolynomials of \textit{non-degenerate exactly-solvable operators}, that is operators of the above type but with the condition $\deg Q_k=k$ for the leading term. We proved that when the degree $n$ of the unique and monic eigenpolynomial $p_n$ tends to infinity, the roots of $p_n$ stay in a compact set in $\mathbb{C}$ and are distributed according to a certain probability measure which is supported by a tree and which depends only on the leading polynomial $Q_k$.\\\\
\textbf{B.} In \cite{B} we studied eigenpolynomials of \textit{degenerate exactly-solvable operators} $(\deg Q_k<k)$. We proved that there exists a unique and monic eigenpolynomial $p_n$ for all sufficiently large values on the degree $n$, and that
 the largest modulus of the roots of $p_n$ tends to infinity when $n\to\infty$. We also presented an explicit conjecture and partial results on the growth of the largest root. Namely,\\\\
\textbf{Conjecture (from \cite{B}).} \textit{Let $T=\sum_{j=1}^{k}Q_jD^j$ be a degenerate exactly-solvable operator of order $k$ and denote by $j_0$ the largest $j$ for which $\deg Q_j=j$. Let $r_n=\max\{|\alpha |:p_n(\alpha)=0\}$, where $p_n$ is the unique and monic $n$th degree eigenpolynomial of $T$. Then 
\begin{displaymath}
\lim_{n\to\infty}\frac{r_n}{n^d}=c_0,
\end{displaymath}
 where $c_0>0$ is a positive constant and} 
\begin{displaymath}
d:=\max_{j\in[j_0+1,k]} \bigg(\frac{j-j_0}{j-\deg Q_j}\bigg).
\end{displaymath}\\
Extensive computer experiments listed in \cite{B} confirm the existence of such a constant $c_0$. Now
consider the scaled eigenpolynomial $q_n(z)=p_n(n^dz)$. We construct the probability measure $\mu_n$ by placing a point mass of size $1/n$ at each zero of $q_n$. Numerical evidence indicates that for each degenerate exactly-solvable operator $T$, the sequence $\{\mu_n\}$ converges weakly to a probability measure $\mu_T$ which is (compactly) supported by a tree. In \cite{B} we deduced the algebraic equation satisfied by the Cauchy transform of $\mu_T$.\footnote{It remains to prove the existence of $\mu_T$ and to describe its support explicitly.} Namely, let $T=\sum_{j=1}^{k}Q_j(z)D^j=\sum_{j=1}^{k}\big(\sum_{i=0}^{\deg Q_j}q_{j,i}z^i\big)D^j$ and denote by $j_0$ the largest $j$ for which $\deg Q_j=j$. Assuming wlog that $Q_{j_0}$ is monic, i.e. $q_{j_0,j_0}=1$, we have 
$$z^{j_0}C^{j_0}(z)+\sum_{j\in A}q_{j,\deg Q_j}z^{\deg Q_j}C^j(z)=1,$$
where $C(z)=\int\frac{d\mu_T(\zeta )}{z-\zeta}$ is the Cauchy transform of $\mu_T$ and $A=\{j:(j-j_0)/(j-\deg Q_j)=d\}$, where $d$ is defined in the conjecture.
Below we present some typical pictures of the roots of the scaled eigenpolynomial $q_n(z)=p_n(n^dz)$.\\\\
\begin{tabular}{ccccc}
\textit{Fig.1:} & & \textit{Fig.2:} & & \textit{Fig.3:}\\
\includegraphics[width=2.8cm]{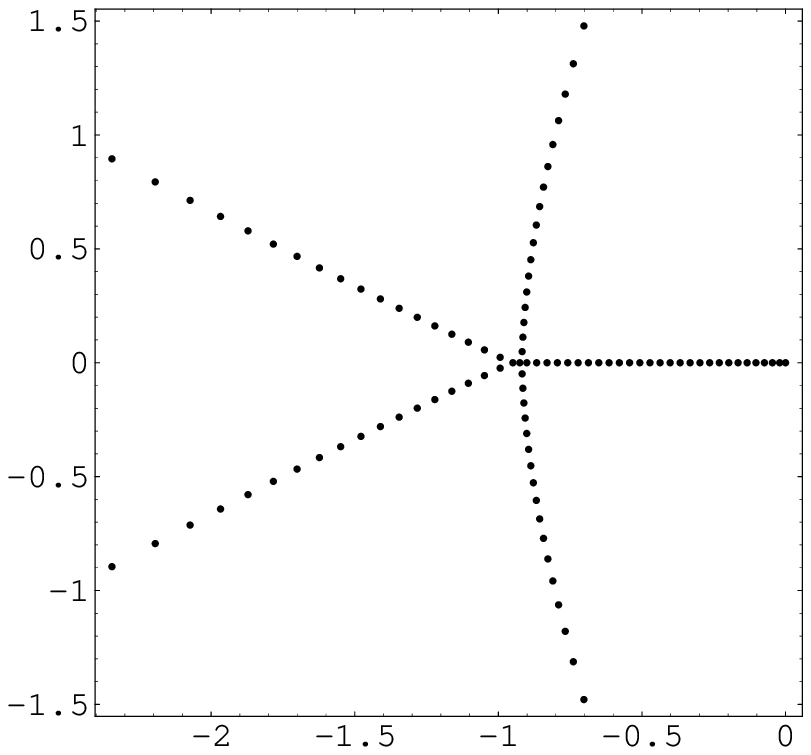} & &
\includegraphics[width=2.8cm]{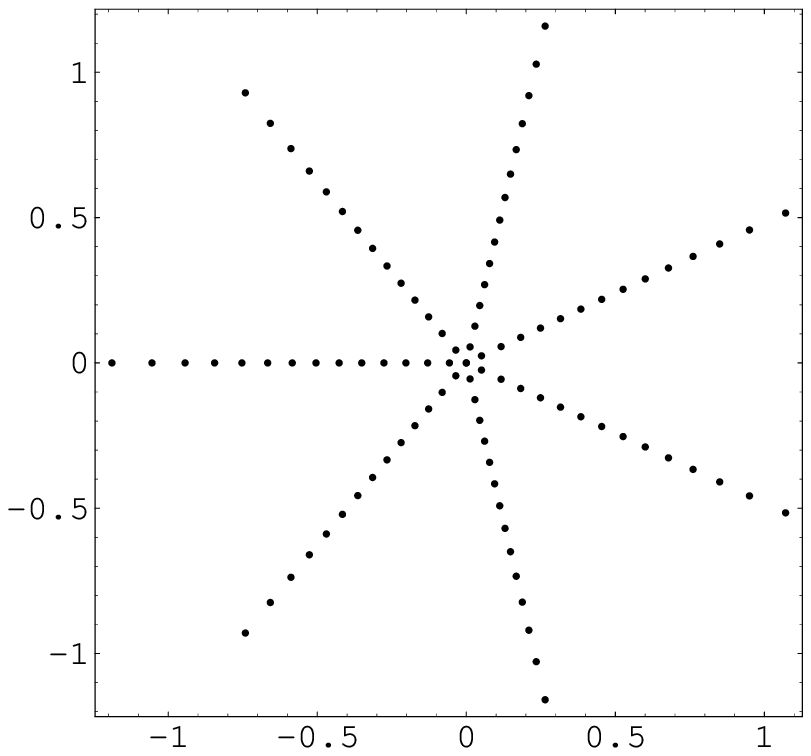} & &
\includegraphics[width=2.8cm]{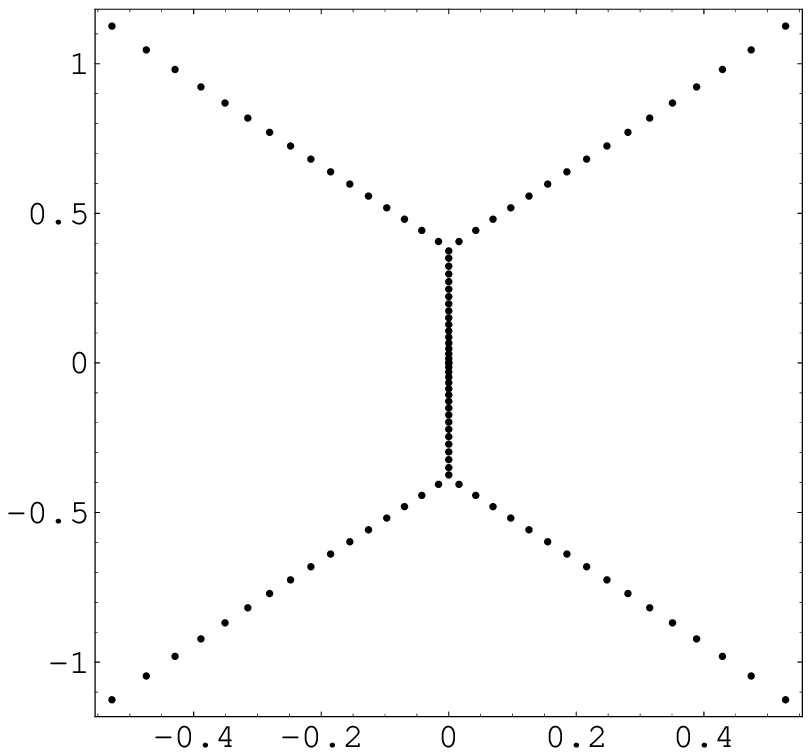}\\
roots of & & roots of & & roots of\\
 $q_{100}(z)=p_{100}(100z)$ & & $q_{100}(z)=p_{100}(100z)$
 & & $q_{100}(z)=p_{100}(100z)$ 
\end{tabular}\\\\\\
\textit{Fig.1:} $T_1=zD+zD^2+zD^3+zD^4+zD^5$.\\
\textit{Fig.2:} $T_2=z^2D^2+D^7$.\\
\textit{Fig.3:} $T_3=z^3D^3+z^2D^4+zD^5$.\\

In this paper we extend the results from \cite{B} by establishing a lower bound for $r_n$ for \textit{all} degenerate exactly-solvable operators and which supports the above conjecture.\footnote{It is still an open problem to prove the upper bound.} This is our main result:\\\\
\textbf{Main Theorem.} \textit{Let $T=\sum_{j=1}^{k}Q_jD^{j}$ be a degenerate exactly-solvable operator and denote by $j_0$ the largest $j$ for which $\deg Q_j=j$. Let $p_n$ be the unique and monic $n$th degree eigenpolynomial of $T$ and 
 $r_n=\max\{|\alpha|:p_n(\alpha)=0\}$. Then there exists a positive constant $c>0$ such that} $$\lim_{n\to\infty}\frac{r_n}{n^d}\geq c,$$ \textit{where}
$$d:=\max_{j\in[j_0+1,k]}\big(\frac{j-j_0}{j-\deg Q_j}\big).$$\\\\

\textbf{Acknowledgements.} The authors are greatly obliged to Professor Boris Shapiro for introducing us to this very fascinating subject. Our research was supported by Stockholm University. 
\section{Proofs} 
\textbf{Lemma 1.} \textit{For any monic polynomial $p(z)$ of degree $n\geq 2$ for which all the zeros are contained in a disc of radius $A\geq 1$,
there exists an integer $n(j)$ and an absolute constant $C_j$ depending only on $j$, such that for every $j\geq 1$ and every $n\geq n(j)$ we have}
\begin{equation}\label{doubleineq}
\frac{1}{C_j}\cdot\frac{n^j}{A^j}\leq\bigg|\bigg|\frac{p^{(j)}(z)}{p(z)}\bigg|\bigg|_{2A}\leq C_j\cdot\frac{n^j}{A^j}
\end{equation}
\textit{where $p^{(j)}(z)$ denotes the $j$th derivative of $p(z)$, and where we have used the maximum norm $||p(z)||_{2A}=\max_{|z|=2A}|p(z)|$.}\\\\
\textbf{Remark.} The right-hand side of the above inequality actually holds for all $n\geq 2$, whereas the left-hand side holds for all $n\geq n(j)$.\\\\
\textbf{Proof.} To obtain the inequality on the \textit{right-hand side} we use the notation 
 $p(z)=\prod_{i=1}^{n}(z-\alpha_i)$ where by assumption $|\alpha_i|\leq A$ for every complex root of $p(z)$. Then $p^{(j)}(z)$ is the sum of $n(n-1)\cdots (n-j+1)$ terms, each being the product of $(n-j)$ factors $(z-\alpha_i)$.\footnote{Differentiating $p(z)=\prod_{i=1}^{n}(z-\alpha_i)$ once yields ${n\choose 1}=n$ terms each term being a product of $(n-1)$ factors $(z-\alpha_i)$, differentiating once again we obtain $n {n-1\choose 1}=n(n-1)$ terms, each being the product of $(n-2)$ factors $(z-\alpha_i)$, etc.} Thus $p^{(j)}(z)/p(z)$ is the sum of $n(n-1)\cdots (n-j+1)$ terms, each equal to $1$ divided by a product consisting of $n-(n-j)=j$ factors $(z-\alpha_i)$. If $|z|=2A$ we get $|z-\alpha_i|\geq A \Rightarrow \frac{1}{|z-\alpha_i|}\leq \frac{1}{A}$, and thus $$\big|\big|\frac{p^{(j)}}{p}\big|\big|_{2A}\leq\frac{n(n-1)\cdots (n-j+1)}{A^j}\leq C_j\cdot\frac{n^j}{A^j}.$$ Here we can choose $C_j=1$ for all $j$, but we refrain from doing this since we will need $C_j$ large enough to obtain the constant $1/C_j$ in the left-hand side inequality. To prove the \textit{left-hand side} inequality we will need inequalities (i)-(iv) below, where we need (i) to prove (ii), and we need (ii) and (iii) to prove (iv), from which the left-hand side inequality of this lemma follows.\\\\
For every $j\geq 1$ we have\\\\
(i)$\qquad\big|\big|\frac{d}{dz}\big(\frac{p^{(j)}(z)}{p(z)}\big)\big|\big|_{2A}\leq j\cdot\frac{n^j}{A^{j+1}}$.\\\\
For every $j\geq 1$ there exists a positive constant $C_j'$ depending only on $j$, such that\\\\
(ii)$\qquad \bigg|\bigg|\frac{p^{(j)}}{p}-\frac{(p')^j}{p^j}\bigg|\bigg|_{2A}\leq C_j'\cdot\frac{n^{j-1}}{A^j}$.\\\\
(iii)$\qquad \big|\big|\frac{p'}{p}\big|\big|_{2A}\geq\frac{n}{3A}$.\\\\
For every $j\geq 1$ there exists a positive constant $C_j''$ and some integer $n(j)$ such that
for all $n\geq n(j)$ we have\\\\
(iv)$\qquad\big|\big|\frac{p^{(j)}}{p}\big|\big|_{2A}\geq C_j''\cdot\frac{n^j}{A^j}$.\\\\
To prove (i), let $p(z)=\prod_{i=1}^{n}(z-\alpha_i)$, where $|\alpha_i|\leq A$ for each complex root $\alpha_i$ of $p(z)$. Then again $p^{(j)}(z)/p(z)$ is the sum of $n(n-1)\cdots (n-j+1)$ terms and each term equals $1$ divided by a product consisting of $j$ factors $(z-\alpha_i)$. Differentiating \textit{each} such term we obtain a sum of $j$ terms each being on the form $(-1)$ divided by a product consisting of $(j+1)$ factors $(z-\alpha_i)$.\footnote{With $D=d/dz$ consider for example $D\frac{1}{\prod_{i=1}^{j}(z-\alpha_i)}=\frac{-1\cdot D\prod_{i=1}^{j}(z-\alpha_i)}{\prod_{i=1}^{j}(z-\alpha_i)^2}$, which is a sum of $j$ terms, each being on the form $(-1)$ divided by a product consisting of $2j-(j-1)=(j+1)$ factors $(z-\alpha_i)$.} 
Thus $\frac{d}{dz}\big(\frac{p^{(j)}(z)}{p(z)}\big)$ is a sum consisting of 
$j\cdot n(n-1)\cdots (n-j+1)$ terms, each on the form $(-1)$ divided by $(j+1)$ factors $(z-\alpha_i)$. Using $\frac{1}{|z-\alpha_i|}\leq \frac{1}{A}$ for $|z|=2A$ since $|\alpha_i|\leq A$ for all $i\in[1,n]$, we thus get 
\begin{displaymath}
\big|\big|\frac{d}{dz}\big(\frac{p^{(j)}(z)}{p(z)}\big)\big|\big|_{2A}\leq 
\frac{j\cdot n(n-1)\cdots (n-j+1)}{A^{j+1}}\leq j\cdot\frac{n^j}{A^{j+1}}.
\end{displaymath}\\\\
To prove (ii) we use (i) and induction over $j$. The case $j=1$ is trivial since 
$\frac{p'}{p}-\frac{(p')^1}{p^1}=0$. If we put $j=1$ in (i) we get 
$\big|\big|\frac{d}{dz}\big(\frac{p'}{p}\big)\big|\big|_{2A}\leq\frac{n}{A^2}$. But $\frac{d}{dz}\big(\frac{p'}{p}\big)=
\frac{p^{(2)}}{p}-\frac{(p')^2}{p^2}$, and thus $\big|\big|\frac{p^{(2)}}{p}-\frac{(p')^2}{p^2}\big|\big|\leq \frac{n}{A^2}$, so (ii) holds for 
$j=2$. We now proceed by induction. Assume that (ii) holds for some $j=p\geq 2$, i.e.
$\big|\big|\frac{p^{(p)}}{p}-\frac{(p')^p}{p^p}\big|\big|_{2A}\leq 
C_p'\cdot\frac{n^{p-1}}{A^p}$. 
Also note that with $j=p$ in (i) we have 
\begin{displaymath}
\big|\big|\frac{p^{(p+1)}}{p}-\frac{p^{(p)}\cdot p'}{p^2}\big|\big|_{2A}=
\big|\big|\frac{d}{dz}\big(\frac{p^{(p)}}{p}\big)\big|\big|_{2A}\leq p\cdot\frac{n^p}{A^{p+1}},
\end{displaymath}
and also $||\frac{p'}{p}||_{2A}\leq\frac{n}{A}$ (from the right-hand side inequality of this lemma).
Thus we have 
\begin{eqnarray*}
\bigg|\bigg|\frac{p^{(p+1)}}{p}-\frac{(p')^{p+1}}{p^{p+1}}\bigg|\bigg|_{2A}&=&
\bigg|\bigg|\frac{p^{(p+1)}}{p}-\frac{p^{(p)}\cdot p'}{p^2}+\frac{p^{(p)}\cdot p'}{p^2}-\frac{(p')^{p+1}}{p^{p+1}}\bigg|\bigg|_{2A}\\
&\leq&
\bigg|\bigg|\frac{p^{(p+1)}}{p}-\frac{p^{(p)}\cdot p'}{p^2}\bigg|\bigg|_{2A}+\bigg|\bigg|\frac{p'}{p}\bigg(\frac{p^{(p)}}{p}-\frac{(p')^p}{p^p}\bigg)\bigg|\bigg|_{2A}\\
&\leq&
p\cdot\frac{n^p}{A^{p+1}}+\frac{n}{A}\cdot C_p'\cdot\frac{n^{p-1}}{A^p}\\
&=&(p+C_p')\cdot\frac{n^p}{A^{p+1}}=
C_{p+1}'\cdot\frac{n^p}{A^{p+1}}.
\end{eqnarray*}\\\\
To prove (iii) observe that $\frac{p'(z)}{p(z)}=\sum_{i=1}^{n}\frac{1}{(z-\alpha_i)}
=\sum_{i=1}^{n}\frac{1}{z}\cdot\frac{1}{1-\frac{\alpha_i}{z}}$. By assumption $|\alpha_i|\leq A$ for all complex roots $\alpha_i$ of $p(z)$, so for $|z|=2A$ we have $|\frac{\alpha_i}{z}|\leq \frac{A}{2A}=\frac{1}{2}$ for all $i\in[1,n]$. Writing $w_i=\frac{1}{1-\frac{\alpha_i}{z}}$ we obtain
\begin{displaymath}
|w_i-1|=\bigg|\frac{1}{1-\frac{\alpha_i}{z}}-\frac{1-\frac{\alpha_i}{z}}{1-\frac{\alpha_i}{z}}\bigg|=
\frac{\big|\frac{\alpha_i}{z}\big|}{\big|1-\frac{\alpha_i}{z}\big|}
\leq \frac{1}{2}|w_i|, 
\end{displaymath}
which implies
\begin{displaymath}
 Re\bigg(\frac{1}{1-\frac{\alpha_i}{z}}\bigg)=Re(w_i)\geq \frac{2}{3}\quad\forall i\in[1,n]\Rightarrow Re\bigg(\sum_{i=1}^{n}\frac{1}{1-\frac{\alpha_i}{z}}\bigg)\geq\frac{2n}{3}.
\end{displaymath}
Thus 
\begin{eqnarray*}
\bigg|\bigg|\frac{p'(z)}{p(z)}\bigg|\bigg|_{2A}&= &
\max_{|z|=2A}\bigg|\frac{p'(z)}{p(z)}\bigg|=\max_{|z|=2A}\frac{1}{|z|}\cdot\bigg|\sum_{i=1}^{n}\frac{1}{1-\frac{\alpha_i}{z}}\bigg|\\
&\geq & \frac{1}{2A}\cdot \bigg|\sum_{i=1}^{n}\frac{1}{1-\frac{\alpha_i}{z}}\bigg|_{2A}\geq
\frac{1}{2A}\cdot Re\bigg(\sum_{i=1}^{n}\frac{1}{1-\frac{\alpha_i}{z}}\bigg)\\
&\geq &\frac{n}{3A}.
\end{eqnarray*}\\\\
To prove (iv) we note that from (iii) we obtain
  $\big|\big|\big(\frac{p'}{p}\big)^j\big|\big|_{2A}\geq\frac{n^j}{3^jA^j}$, and this together with (ii) yields
\begin{eqnarray*}
\bigg|\bigg|\frac{p^{(j)}}{p}\bigg|\bigg|_{2A}&=&\bigg|\bigg|\bigg(\frac{p'}{p}\bigg)^j+\frac{p^{(j)}}{p}
-\bigg(\frac{p'}{p}\bigg)^j\bigg|\bigg|_{2A}\geq\bigg|\bigg|\bigg(\frac{p'}{p}\bigg)^j\bigg|\bigg|_{2A}-
\bigg|\bigg|\frac{p^{(j)}}{p}-\bigg(\frac{p'}{p}\bigg)^j\bigg|\bigg|_{2A}\\
&\geq &\frac{n^j}{3^jA^j}-C_j'\cdot\frac{n^{j-1}}{A^j}=\frac{n^j}{A^j}\bigg(\frac{1}{3^j}-\frac{C_j'}{n}\bigg)\geq
C_j''\cdot\frac{n^j}{A^j},
\end{eqnarray*}
where $C_j''$ is a positive constant such that $C_j''\leq\big(\frac{1}{3^j}-\frac{C_j'}{n}\big)$ for all $n\geq n(j)$.\\\\
The left-hand side inequality in this lemma now follows from (iv) if we choose the constant $C_j$ on right-hand side inequality so large that $\frac{1}{C_j}\leq C_j''$.$\hfill\square$\\\\
To prove Main Theorem we will need the following lemma, which follows from Lemma 1:\\\\
\textbf{Lemma 2.} Let $0<s<1$ and $d>0$ be real numbers. Let $p(z)$ be any monic polynomial of degree $n\geq 2$ such that all its zeros are contained in a disc of radius $A=s\cdot n^d$, and let $Q_j(z)$ be arbitrary polynomials. Then there exists some positive integer $n_0$ and positive constants $K_j$ such that
\begin{displaymath}
\frac{1}{K_j}\cdot n^{d(\deg Q_j-j)+j}\cdot\frac{s^{\deg Q_j}}{s^j}
\leq\bigg|\bigg|Q_j(z)\cdot\frac{p^{(j)}}{p}\bigg|\bigg|_{2sn^{d}}
\leq K_j\cdot n^{d(\deg Q_j-j)+j}\cdot\frac{s^{\deg Q_j}}{s^j}
\end{displaymath}
for every $j\geq 1$ and all $n\geq\max(n_0,n(j))$, where $n(j)$ is as in Lemma 1.\\\\
\textbf{Proof.} Let $Q_j(z)=\sum_{i=0}^{\deg Q_j}q_{j,i}z^i$. Then for $|z|=2A>>1$ we have
\begin{displaymath}
|Q(z)|_{2A}=|q_{j,\deg Q_j}|2^{\deg Q_j}A^{\deg Q_j}\bigg(1+O(\frac{1}{A})\bigg).\end{displaymath}
Since $A=s\cdot n^d$ there exists some integer $n_0$ such that $n\geq n_0\Rightarrow A\geq A_0>>1$, and thus by Lemma 1 there exists a positive constant $K_j$ such that the following inequality holds for all $n\geq\max(n(j),n_0)$ and all $j\geq 1$:
\begin{displaymath}
\frac{1}{K_j}\cdot\frac{n^j}{A^j}\cdot A^{\deg Q_j}\leq\bigg|\bigg|Q_j(z)\cdot\frac{p^{(j)}}{p}\bigg|\bigg|_{2A}\leq K_j\cdot\frac{n^j}{A^j}\cdot
A^{\deg Q_j}.
\end{displaymath}
Inserting $A=s\cdot n^d$ in this inequality  
we obtain\\\\
\begin{displaymath}
\frac{1}{K_j}\cdot\frac{n^j}{s^jn^{dj}}\cdot s^{\deg Q_j}n^{d\cdot\deg Q_j}\leq\bigg|\bigg|Q_j(z)\cdot\frac{p^{(j)}}{p}\bigg|\bigg|_{2sn^{d}}
\leq K_j\cdot\frac{n^j}{s^jn^{dj}}\cdot s^{\deg Q_j}n^{d\cdot\deg Q_j}
\end{displaymath}
\begin{displaymath}
\Leftrightarrow
\end{displaymath}
\begin{displaymath}
\frac{1}{K_j}\cdot n^{d(\deg Q_j-j)+j}\cdot\frac{s^{\deg Q_j}}{s^j}
\leq\bigg|\bigg|Q_j(z)\cdot\frac{p^{(j)}}{p}\bigg|\bigg|_{2sn^{d}}
\leq K_j\cdot n^{d(\deg Q_j-j)+j}\cdot\frac{s^{\deg Q_j}}{s^j}
\end{displaymath}
for every $j\geq 1$ and all $n\geq\max(n_0,n(j))$.$\hfill\square$\\\\
\textbf{Proof of Main Theorem.}  Let $d=\max_{j\in[j_0+1,k]}\big(\frac{j-j_0}{j-\deg Q_j}\big)$ where $j_0$ is the largest $j$ for which $\deg Q_j=j$ in the degenerate exactly-solvable operator $T=\sum_{j=1}^{k}Q_jD^j$, where $Q_j(z)=\sum_{i=0}^{\deg Q_j}q_{j,i}z^i$. Let $p_n(z)$ be the $n$th degree unique and monic eigenpolynomial of $T$ and denote by 
$\lambda_n$ the corresponding eigenvalue. 
Then the eigenvalue equation can be written
\begin{equation}\label{eigenvalueeq}
\sum_{j=1}^{k}Q_j(z)\cdot\frac{p_n^{(j)}(z)}{p_n(z)}=\lambda_n
\end{equation}
where $\lambda_n=\sum_{j=1}^{j_0}q_{j,j}\cdot\frac{n!}{(n-j)!}$.
We will now use the result in Lemma 2 to estimate each term in (\ref{eigenvalueeq}).\\

* Denote by $j_m$ the largest $j$ for which $d$ is attained. Then $d=(j_m-j_0)/(j_m-\deg Q_{j_m})\Rightarrow d(\deg Q_{j_m}-j_m)+j_m=j_0$, and $j_m-\deg Q_{j_m}=(j_m-j_0)/d$. By Lemma 2 we have:
\begin{equation}
\frac{1}{K_{j_m}}\cdot n^{j_0}\cdot\frac{1}{s^{\frac{j_m-j_0}{d}}}
\leq\bigg|\bigg|Q_{j_m}(z)\cdot\frac{p^{(j_m)}}{p}\bigg|\bigg|_{2sn^{d}}
\leq K_{j_m}\cdot n^{j_0}\cdot\frac{1}{s^{\frac{j_m-j_0}{d}}}.
\end{equation}
Note that the exponent of $s$ is positive since $j_m>j_0$ and $d>0$. In what follows we will only need the left-hand side of the above inequality. \\

* Consider the remaining (if there are any) $j_0<j<j_m$ for which $d$ is attained. For such $j$ we have (using the right-hand side inequality of Lemma 2):
\begin{eqnarray}
\bigg|\bigg|Q_{j}(z)\cdot\frac{p^{(j)}}{p}\bigg|\bigg|_{2sn^{d}}
&\leq& K_{j}n^{j_0}\cdot\frac{1}{s^{\frac{j-j_0}{d}}}= K_{j}n^{j_0}\cdot\frac{1}{s^{\frac{j_m-j_0}{d}}}\cdot
s^{\frac{j_m-j}{d}}\nonumber\\
&\leq& K_jn^{j_0}\cdot\frac{1}{s^{\frac{j_m-j_0}{d}}}\cdot s^{1/d}
\end{eqnarray}
where we have used that $(j_m-j)\geq 1$ and $s<1\Rightarrow s^{(j_m-j)/d}\leq s^{1/d}$.\\

* Consider all $j_0<j\leq k$ for which $d$ is \textit{not} attained. Then $(j-\deg Q_j)>0$ and $(j-j_0)/(j-\deg Q_j)<d\Rightarrow 
d(\deg Q_j-j)+j<j_0$ and we can write $d(\deg Q_j-j)+j\leq j_0-\delta$ where $\delta>0$. Then we have:
\begin{eqnarray}
\bigg|\bigg|Q_j(z)\cdot\frac{p^{(j)}}{p}\bigg|\bigg|_{2sn^{d}}
&\leq& K_j\cdot n^{d(\deg Q_j-j)+j}\cdot\frac{s^{\deg Q_j}}{s^j}\leq K_j\cdot n^{j_0-\delta}\cdot\frac{s^{\deg Q_j}}{s^j}\nonumber\\
&\leq& K_j\cdot n^{j_0-\delta}\cdot\frac{1}{s^k}, 
\end{eqnarray}
where the last inequality follows since $\deg Q_j\geq 0\Rightarrow s^{\deg Q_j}\leq s^0=1$ and $j\leq k\Rightarrow s^j\geq s^k$ since $0<s<1$.\\

* For $j=j_0$ by definition $\deg Q_{j_0}=j_0$ and thus:
\begin{equation}
\bigg|\bigg|Q_{j_0}(z)\cdot\frac{p^{(j_0)}}{p}\bigg|\bigg|_{2sn^{d}}
\leq K_{j_0}\cdot n^{d(\deg Q_{j_0}-j_0)+j_0}\cdot\frac{s^{\deg Q_{j_0}}}{s^{j_0}}=K_{j_0}\cdot n^{j_0}.
\end{equation}
\\

* Now consider all $1\leq j\leq j_0-1$. Since $n\geq n_0\Rightarrow A=sn^d>>1$ we get
$(sn^d)^{j-\deg Q_j}\geq 1$ and thus:
\begin{eqnarray}
\bigg|\bigg|Q_j(z)\cdot\frac{p^{(j)}}{p}\bigg|\bigg|_{2sn^{d}}
&\leq &K_j\cdot n^{d(\deg Q_j-j)+j}\cdot\frac{s^{\deg Q_j}}{s^j}=K_j\cdot n^j\cdot (sn^d)^{(\deg Q_j-j)}\nonumber\\
&=&K_j\cdot n^j\cdot\frac{1}{(sn^d)^{j-\deg Q_j}}\leq 
K_j\cdot n^j\leq K_j\cdot n^{j_0-1}.
\end{eqnarray}\\

* Finally we estimate the eigenvalue $\lambda_n=\sum_{i=1}^{j_0}q_{j,j}\cdot\frac{n!}{(n-j)!}$, which grows as $n^{j_0}$ for large $n$, since there
exists an integer $n_{j_0}$ and some positive constant $K_{j_0}'$ such that for all $n\geq n_{j_0}$ we obtain: 
\begin{eqnarray}
|\lambda_n|&\leq &\sum_{j=1}^{j_0}|q_{j,j}|\cdot\frac{n!}{(n-j)!}=|q_{j_0,j_0}|\cdot\frac{n!}{(n-j_0)!}\bigg[1+
\sum_{1\leq j<j_0}\bigg|\frac{q_{j,j}}{q_{j_0,j_0}}\bigg|\cdot\frac{(n-j_0)!}{(n-j)!}\bigg]\nonumber\\
&\leq & K_{j_0}'\cdot n^{j_0}.
\end{eqnarray}
Finally we rewrite the eigenvalue equation (\ref{eigenvalueeq}) as follows:
\begin{displaymath}
Q_{j_m}(z)\cdot\frac{p_n^{(j_m)}(z)}{p_n(z)}=\lambda_n+\sum_{j\neq j_m}Q_j(z)\frac{p_n^{(j)}(z)}{p_n(z)}.
\end{displaymath}
Applying inequalities (5)-(9) to this we obtain
\begin{eqnarray}\label{RHS}
\bigg|\bigg|Q_{j_m}\cdot\frac{p_n^{(j_m)}(z)}{p_n(z)}\bigg|\bigg|_{2sn^d}&\leq & |\lambda_n|+\sum_{j\neq j_m}\bigg|\bigg|Q_j\frac{p_n^{(j)}(z)}{p_n(z)}\bigg|\bigg|_{2sn^d}\nonumber\\
&\leq & K_{j_0}'n^{j_0}+K_{j_0}n^{j_0}+\sum_{1\leq j<j_0}K_jn^{j_0-1}\nonumber\\
&+& \sum_{j_0<j\leq k:\atop (\frac{j-j_0}{j-\deg Q_j}<d}) K_j\frac{n^{j_0-\delta}}{s^k}+
\sum_{j_0<j<j_m\atop (\frac{j-j_0}{j-\deg Q_j}=d}) K_jn^{j_0}\frac{s^{1/d}}{s^{\frac{j_m-j_0}{d}}}\nonumber\\
 &\leq & K\cdot n^{j_0}+
K\cdot\frac{n^{j_0-\delta}}{s^k}+
 K\cdot n^{j_0}\frac{s^{1/d}}{s^{\frac{j_m-j_0}{d}}}
\end{eqnarray}
 for all $n\geq\max(n_0, n(j), n_{j_0})$, where 
 $K$ is some positive constant
and $0<s<1$. For the term on the left-hand side of the rewritten eigenvalue equation above we obtain using (4) the following estimation:
\begin{eqnarray}\label{LHS}
 \frac{1}{K}\cdot n^{j_0}\cdot\frac{1}{s^{\frac{j_m-j_0}{d}}}   
\leq
 \frac{1}{K_{j_m}}\cdot n^{j_0}\cdot\frac{1}{s^{\frac{j_m-j_0}{d}}}   
\leq\bigg|\bigg|Q_{j_m}\cdot\frac{p_n^{(j_m)}(z)}{p_n(z)}\bigg|\bigg|_{2sn^d}
\end{eqnarray}
for some constant $K\geq K_{j_m}$ which also satisfies (\ref{RHS}). 
Now combining (\ref{RHS}) and (\ref{LHS}) we get
\begin{equation*}
\frac{1}{K}\cdot n^{j_0}\cdot\frac{1}{s^{\frac{j_m-j_0}{d}}}\leq K\cdot n^{j_0}+
K\cdot\frac{n^{j_0-\delta}}{s^k}+ K\cdot n^{j_0}\frac{s^{1/d}}{s^{\frac{j_m-j_0}{d}}}.
\end{equation*}
Dividing this inequality by $n^{j_0}$ and multiplying by $K$ we have
\begin{equation*}
\frac{1}{s^{\frac{j_m-j_0}{d}}}\leq K^2+
K^2\cdot\frac{1}{n^{\delta}}\cdot\frac{1}{s^k}+ K^2\cdot\frac{s^{1/d}}{s^{\frac{j_m-j_0}{d}}}.
\end{equation*}
\begin{displaymath}
\Leftrightarrow
\end{displaymath}
\begin{displaymath}
\frac{1}{s^w}\leq K^2+\frac{K^2}{s^k}\cdot\frac{1}{n^{\delta}}+K^2\cdot\frac{s^{1/d}}{s^w}
\end{displaymath}
\begin{displaymath}
\Leftrightarrow
\end{displaymath}
\begin{equation}\label{ineq}
\frac{1}{s^w}[1-K^2\cdot s^{1/d}]\leq K^2+\frac{K^2}{s^k}\cdot\frac{1}{n^{\delta}}.
\end{equation}
where $w=(j_m-j_0)/d>0$.\\

 In what follows we will obtain a contradiction to this inequality for some small properly chosen $0<s<1$ and all sufficiently large $n$. 
Since $j_m\in[j_0+1,k]$ we have $w=(j_m-j_0)/d\geq 1/d$, and since $s<1$ we get $s^w\leq s^{1/d}\Rightarrow
 1/s^w\geq 1/s^{1/d}$. \textbf{Now choose $s^{1/d}=\frac{1}{4K^2}$, where $K$ is the constant in (\ref{ineq})}. Then 
estimating the left-hand side of (\ref{ineq}) we get
\begin{eqnarray*}
\frac{1}{s^w}[1-K^2\cdot s^{1/d}]\geq \frac{1}{s^{1/d}}[1-K^2\cdot s^{1/d}]=4K^2-K^2=3K^2
\end{eqnarray*}
and thus from  (\ref{ineq}) we have 
\begin{displaymath}
3K^2\leq \frac{1}{s^w}[1-K^2\cdot s^{1/d}]\leq K^2+\frac{K^2}{s^k}\cdot\frac{1}{n^{\delta}}
\end{displaymath}
\begin{displaymath}
\Leftrightarrow
\end{displaymath}
\begin{displaymath}
2K^2\leq\frac{K^2}{s^k}\cdot\frac{1}{n^{\delta}}
\end{displaymath}
\begin{displaymath}
\Leftrightarrow
\end{displaymath}
\begin{displaymath}
n^{\delta}\leq\frac{1}{2}\cdot\frac{1}{s^k}=\frac{1}{2}(2K)^{2dk}. 
\end{displaymath}
We therefore obtain a contradiction to this inequality, and hence to inequality (\ref{ineq}) and thus to the eigenvalue equation (\ref{eigenvalueeq}),
if $n^{\delta}>\frac{1}{2}(2K)^{2dk}$ and $s=1/(2K)^{2d}$, and consequently all roots of $p_n$ cannot be contained in a disc of radius $s\cdot n^d$ for such choices on $s$ and $n$, whence $r_n>s\cdot n^d$ where $r_n$ denotes the largest modulus of all roots of $p_n$, so clearly there exists some positive constant $c$ such that 
$\lim_{n\to\infty}\frac{r_n}{n^d}\geq c$.$\hfill\square$
\section{Open Problems and Conjectures}
\subsection{The upper bound}
Based upon numerical evidence from computer experiments (some of which is presented in \cite{B}) we are led to assert that there
exists a constant $C_0$, which depends on
$T$ only, such that
\begin{equation}\label{upperbound}
r_n\leq C_0\cdot n^d
\end{equation}
holds for all sufficiently large integers $n$. We refer to this as the \textbf{upper-bound conjecture}. Computer experiments confirm that the zeros of the scaled eigenpolynomial $q_n(z)=p_n(n^dz)$ are contained in a compact set when $n\to\infty$. 
\subsection{The measures $\{\mu_n\}$}
Consider the sequence of discrete
probability measures
$$
\mu_n=\frac{1}{n}\sum_{\nu=1}^{\nu=n}
\delta(\frac{\alpha_\nu}{n^d})
$$
where $\alpha_1,\ldots,\alpha_n$ are the roots of the
eigenpolynomial
$p_n$ and $d$ is as in Definition 1.
Assuming \eqref{upperbound} the supports of
$\{\mu_n\}$
stay in a compact set in $\mathbb{C}$.
Next, by a \textbf{tree} we mean a connected
compact subset
$\Gamma$ of $\mathbb{C}$ which consists of a finite union of
real-analytic curves and where
$\hat{\mathbb{C}}\setminus\Gamma$
is simply connected (here
$\hat{\mathbb{C}}=\mathbb{C}\cup \infty$ is the extended complex
plane). Computer experiments from \cite{B}
lead us to the following\\\\
\textbf{Conjecture 1.}
For each operator
$T$ the sequence $\{\mu_n\}$
converges weakly to a probability measure
$\mu_T$ which is supported on a certain tree
$\Gamma_T$.
\subsection{The determination of
$\mu_T$}
Given  $T=\sum_{j=1}^{k}Q_j(z)D^j$ and $Q_j(z)=\sum_{i=0}^{\deg Q_j}q_{j,i}z^i$
we obtain an algebraic function $y_T(z)$
which satisfies the following algebraic equation (also see \cite{B}):
$$
q_{j_0,j_0}\cdot z^{j_0}\cdot y^{j_0}_T(z)+\sum_{j\in J}
q_{j,\deg Q_j}\cdot z^{\deg Q_j}\cdot y^j_T(z)=q_{j_0,j_0},
$$
where $J=\{j:(j-j_0)/(j-\deg Q_j)=d\}$, i.e. the sum is taken over all
$j$ for which  $d$ is attained. 
In addition $y_T$
is chosen to be the unique single-valued branch
which has an expansion
$$
y_T(z)=\frac{1}{z}+ \frac{c_2}{z^2}+\frac{c_3}{z^3}+
\ldots
$$
at $\infty\in\hat{\mathbb{C}}$. 
Let $\mathbb{D}_T$ be the discriminant locus of $y_T$, i.e. this is
a finite set in $\mathbb{C}$
such that the single-valued branch
of $y_T$ in an exterior disc
$|z|>R$ can be continued to
an
(in general multi-valued) analytic function in
$\hat{\mathbb{C}}\setminus\mathbb{D}_T$.
If $\Gamma_T$
is a tree which contains
$\mathbb{D}_T$, we obtain a single-valued branch of
$y_T$ in the simply connected
set
$\Omega_{\Gamma_T}=\hat{\mathbb{C}}\setminus\Gamma_T$.
It is easily seen that this holomorphic function in
$\Omega_{\Gamma_T}$
defines a locally integrable function  in the sense of Lebesgue outside the
nullset $\Gamma_T$.
A tree $\Gamma_T$ which contains
$\mathbb{D}_T$ is called $T$-positive if
the distribution defined by
$$
\nu_{\Gamma_T}=\frac{1}{\pi}\cdot
{\bar\partial y_T/\bar\partial\bar z}
$$
is a probability measure.
\subsection{Main conjecture}
Now we  announce the following
which  is
experimentally confirmed
in \cite{B}:\\\\
\textit{For each operator
$T$, the limiting measure $\mu_T$
in Conjecture 1 exists.  Moreover, its
support
is a $T$-positive
tree $\Gamma_T$
and one has the equality
$\mu_T=\nu_{\Gamma_T}$  which means that when $z\in \hat{\mathbf{C}}\setminus\Gamma_T$ the following holds:}
$$
y_T(z)=\int_{\Gamma_T}\,{d\mu_T(\zeta)\over z-\zeta}.
$$
\textbf{Remark.}
For
{\it{non-degenerate exactly-solvable operators}}
(i.e. when $\deg Q_k=k$)
 it was proved in
\cite{BR} that the roots of all eigenpolynomials
stay in a compact set of $\mathbb{C}$, and the
unscaled sequence of
probability measures
$\{\mu_n\}$ converge to a measure supported on a tree, i.e.
the analogue of the main conjecture above  holds in
the non-degenerate case.


\newpage

\begin{thebibliography}{99}
\bibitem{B} T.~Bergkvist: On Asymptotics of Polynomial Eigenfunctions for Exactly-Solvable Differential Operators, math.SP/0701143, 
to appear in \textit{J. Approx. Th.}. 

\bibitem{BR} T.~Bergkvist and H.~Rullg\aa rd:
On polynomial eigenfunctions for a class of  differential
operators, \textit{Math. Research Letters} \textbf{9}, 153 -- 171 (2002).

\bibitem{BRS} T.~Bergkvist, H.~Rullg\aa rd and B.~Shapiro:
On Bochner-Krall Orthogonal Polynomial Systems,
  \textit{Math.Scand} \textbf{94}, no. 1, 148-154 (2004).

\bibitem{TB} T.~Bergkvist: On generalized Laguerre Polynomials with Real and Complex Parameter, \textit{Research Reports in Mathematics, Stockholm University} No. 2 (2003), available at http://www.math.su.se/reports/2003/2/.

\bibitem{BBS} J.~Borcea, R.~B\o gvad, B.~Shapiro: 
On Rational Approximation of Algebraic Functions, to appear in \textit{Adv. Math}, math. CA /0409353.


\bibitem{MS} G.~M\'asson and B.~Shapiro:
A note on  polynomial eigenfunctions of a hypergeometric type
operator,
\textit{Experimental Mathematics}, \textbf{10}, 609--618.

\bibitem{FGZ} A.~Martinez-Finkelshtein, P.~Martinez-Gonzalez, A.~Zarzo: 
WKB approach to zero distribution of solutions of linear second order differential equations, \textit{J. Comp. Appl. Math.} \textbf{145} (2002), 167-182.

\bibitem{MMO} A.~Martinez-Finkelshtein, P.~Martinez-Gonzalez, R.~Orive: 
On asymptotic zero distribution of Laguerre and generalized Bessel polynomials with varying parameters. Proceedings of the Fifth International Symposium on Orthogonal Polynomials, Special Functions and their Applications (Patras 1999), 
\textit{J. Comput. Appl. Math.} \textbf{133} (2001), no. 1-2, p. 477-487.

\bibitem{AT1} A.~Turbiner: Lie-Algebras and Linear Operators with Invariant Subspaces, \textit{Lie Algebras, Cohomologies and New Findings in Quantum Mechanics} AMS Contemporary Mathematics' series, N. Kamran and P. Olver (Eds.), vol \textbf{160}, 263-310 (1994).

\bibitem{AT2} A.~Turbiner: On Polynomial Solutions of differential equations, 
\textit{J. Math. Phys.} \textbf{33} (1992) p.3989-3994.

\bibitem{AT3} A.~Turbiner: Lie algebras and polynomials in one variable, 
\textit{J. Phys. A: Math. Gen.} \textbf{25} (1992) L1087-L1093.
\end{thebibliography}
\end{document}